\documentclass{amsart}
\usepackage{graphicx}
\usepackage{amsmath,amsthm,amsfonts,amssymb}
\newcommand{\Q}{\mathbb{Q}}
\newcommand{\Z}{\mathbb{Z}}
\newcommand{\A}{\mathbb{A}}
\newcommand{\OK}{\mathcal{O}}
\newcommand{\Gal}{\text{Gal}}

\newcommand{\cp}{\mathcal{C}^\prime}

\newcommand{\Jac}{\text{Jac}(\mathcal{C})}
\newcommand{\Pico}{\text{Pic}^{0}(\mathcal{C})}
\newcommand{\rk}{\text{rk}}
\newcommand{\J}{\mathcal{J}}

\title{On Ranks of Jacobian Varieties in Prime Degree Extensions}
\author{Dave Mendes da Costa}
\email{madjmdc@bris.ac.uk}
\address{{School of Mathematics, University Walk} \\
{Bristol, BS8 1TW, United Kingdom}}%
\date{\today}
\begin{document}

\begin{abstract}
In \protect\cite{TD} it is shown that given an elliptic curve $E$ defined over a number field $K$ then there are infinitely many degree 3 extensions $L/K$ for which the rank of $E(L)$ is larger than $E(K)$. In the present paper we show that the same is true if we replace 3 by any prime number. This result follows from a more general result establishing a similar property for the Jacobian varieties associated with curves defined by an equation of the shape $g(y) = f(x)$ where $f$ and $g$ are polynomials of coprime degree. 
\end{abstract}
\maketitle

\section{Introduction and Statement of Results}

Let $\mathcal{C}$ be a curve (which will be smooth, irreducible and projective unless otherwise stated) defined over a number field $K$ and with genus $g(\mathcal{C}) \geq 1$. We can associate with $\mathcal{C}$ an abelian variety $\Jac$ called the $\textit{Jacobian}$ of $\mathcal{C}$. This variety has dimension equal to $g(\mathcal{C})$ and if $\mathcal{C}(K) \neq \emptyset$ then there is an embedding of $\mathcal{C}$ into $\Pico$ defined over $K$ and an isomorphism, also defined over $K$, between $\Jac$ and $\text{Pic}^{0}(\mathcal{C})$. Let us assume that $\mathcal{C}(K) \neq \emptyset$ so that we can identify $\Jac$ and $\Pico$ throughout. \\
\\
The celebrated Mordell-Weil Theorem tells us that $\Jac(K)$ has the structure of a finitely generated abelian group. We define the $\textit{rank}$ of $\Jac$ to be the number of copies of $\Z$ appearing in $\Jac(K)$ and denote this by $\rk(\mathcal{C} / K)$. In this paper we shall be interested in how $\rk(\mathcal{C} / L)$ behaves as we vary the field $L$. When $\mathcal{C}$ is an elliptic curve and $K = \Q$ then a conjecture of Goldfeld \protect\cite{G} asserts that as we let $L$ range across all quadratic extensions of $\Q$ then the rank should remain the same as $\rk(\mathcal{C} /\Q)$ $50 \%$ of the time and increase by one $50 \%$ of the time with the remaining $0 \%$ accounting for other behaviour. This is as yet unproved however it is known that the rank both increases infinitely often and remains the same infinitely often. From this position it is a natural question to ask if this behaviour persists when we make two generalisations:

\begin{enumerate}
	\item Replacing $\Q$ by $K$ and;
	\item Considering $L/K$ such that $[L:K] = p$ for some prime $p$. \\
\end{enumerate}
In this first case the analogue of Goldfeld's conjecture clashes with other standard conjectures which predict that there are elliptic curves defined over number fields whose rank increases in every quadratic extension (for details, see \protect\cite{DD}). With this in mind we ask two questions:

\newtheorem{Question}{Question}
\begin{Question}
\label{q:1}
Given a curve $\mathcal{C}$, of positive genus, defined over a number field $K$ such that $\mathcal{C}(K) \neq \emptyset$ and a prime $p$, then are there infinitely many extensions $L/K$ with $[L:K]=p$ such that $\rk (\mathcal{C} / L) > \rk (\mathcal{C} / K)$ ? \\
\end{Question}

We will also be interested in the following related question:

\begin{Question}
\label{q:2}
Given a curve $\mathcal{C}$, of positive genus, defined over a number field $K$ such that $\mathcal{C}(K) \neq \emptyset$, is there an $N>0$ such that if $p$ is a prime and $p \geq N$ then are there infinitely many extensions $L/K$ with $[L:K]=p$ such that $\rk (\mathcal{C} / L) > \rk (\mathcal{C} / K)$ ? \\
\end{Question}
\newtheorem{thm}{Theorem}
\newtheorem{cor}{Corollary}

Our main result is a partial answer to Question \protect ~\ref{q:2} for a certain family of curves.

\begin{thm}
\label{th:jac}
Let $\mathcal{C}$ be a smooth, irreducible curve of positive genus defined over a number field $K$ and with $\mathcal{C}(K) \neq \emptyset$. Suppose further that $\mathcal{C}$ is birational to a plane curve $\cp$ of the form

\[\cp : g(y) = f(x) \]
\\
where the degrees of $f$ and $g$ are coprime. Then there is an integer $N(\mathcal{C})>0$, effective and depending on $\mathcal{C}$, such that for all primes $p \geq N(\mathcal{C})$ we have an affirmative answer for Question \protect ~\ref{q:1}. 
\end{thm}

This theorem has some interesting corollaries. First we note that the family includes all hyperelliptic curves of odd degree.

\begin{cor}
\label{th:hyp}
Let $\mathcal{C}$ be birational to the plane curve cut out by the equation 

\[ \cp: y^2 = g(x)\]
\\
where $k=\text{deg}(g)$ is odd. Then there is an $N(\mathcal{C})$ such that for all primes $p \geq N(\mathcal{C})$ Question \protect ~\ref{q:1} has an affirmative answer. What is more, we can take $N(\mathcal{C}) = k+1$.  
\end{cor}

This has the following pleasing corollary.

\begin{cor}
\label{th:ell}
Question 1 has an affirmative answer for every prime $p$ when $\mathcal{C}$ is an elliptic curve.
\end{cor}

$\textbf{Acknowlegdements.}$ I would like to thank Tim Dokchister for his encouragement and his reading of earlier drafts. I would also like to express my gratitude to Vladamir Dokchitser and Dan Loughran for answering several questions and to Ariyan Javanpeykar for reminding me that there is more in mathematics than elliptic curves! Finally I would also like to thank the Ecole Normale Superieure in Paris for providing me with a place to work for the last year. This work was supported by an EPSRC Doctoral Training Award.

\section{Strategy of proof}

The strategy we shall employ for proving these theorems is wholly inspired by the paper \protect\cite{TD} of Tim Dokchitser where he proves that Question \ref{q:1} is true for elliptic curves over number fields when the prime $p$ is 3. Moreover, he shows that this is true even if one restricts the number fields $L$ to be of the form $K(\sqrt[3]{m})$ for some $m \in K$.\\
\\
The idea is as follows. Let us call an element of $L \backslash K$ a $\textit{strictly-L}$ element of $L$. We note that for a prime $p$, the truth of Question \ref{q:1} is equivalent to there being infinitely many extensions $L/K$ of degree $p$ such that $\Jac(L)$ contains a strictly-$L$ point, i.e.,  $\Jac(L) \backslash \Jac(K) \neq \emptyset$. This is essentially shown in \protect\cite{TD} but we shall prove it now for the sake of completeness.

\newtheorem{lem}{Lemma}

\begin{lem}
\label{lem:1}
Let $\mathcal{C}$ be a smooth, projective curve defined over $K$ with $g(\mathcal{C}) \geq 1$ and $\text{p}$ a prime number. Let $\J = \Jac$. Then there are infinitely many degree $\text{p}$ extensions $L/K$ such that $rk (\mathcal{C} / L) > \rk(\mathcal{C} / K)$ if and only if there are infinitely many such $L/K$ such that $\J(L)$ has a strictly-L point.  
\end{lem}

\proof{
Clearly if the rank increases then a new point has been obtained, so one direction is clear. For the other direction we note that the only way in which a new point does not lead to an increase in rank is if the point divides a point in $\J(K)$. We claim this can only happen in finitely many degree $p$ extensions. First of all, let $F$ be the compositum of all the degree $p$ fields. Then the torsion of $\J(F)$ is finite since the residue field for each prime of $F$ is finite. Hence there are only finitely many degree $p$ extensions in which we obtain new torsion.\\
\\
Having dealt with new points which divide the identity we need to consider points which divide other points in $E(K)$. Such points in an extension $L/K$ would lead to the map

\[ f : \frac{\J(K)}{\ell \J(K)} \rightarrow \frac{\J(M)}{\ell \J(M)} \]
\\
failing to be injective for some prime $\ell$ where $M$ is the Galois closure of $L$. The kernel of $f$ is contained in the cohomology group $H^1(\Gal (M/K), \J[\ell])$. Since $\Gal (M/L)$ has order dividing $p!$ we see that if $\ell > p$ then this cohomology group vanishes implying that $f$ is injective. For $\ell \leq p$ it suffices to observe that there are only finitely many degree $p$ extensions $L$ in which we can gain a point $Q \in \J(L)$ such that $\ell Q = R$ for some $R \in E(K)$. To see this note that $\ell Q = R = a_1 P_1 + \ldots + a_r P_r$ where the $P_i$ generate $\J(K)$ and by repeatedly subtracting multiples of $\ell P_i$ for each $i$ we can assume that $a_i < \ell$. There are only finitely many such $R$ and hence only finitely many degree $p$ extensions in which they become divisible by $\ell$.\\
\\
Thus we see that an infinitude of new points all in different degree $p$ extensions implies that the rank must increase in infinitely many of those extensions.
\qed     \\    
}
\\
The next step is to find infinitely many $L/K$ of degree $p$ such that $\J(L)$ has a strictly-$L$ point. This is achieved by constructing such points on $\mathcal{C}$ itself and then carrying them over to $\J$ via the embedding 

\begin{eqnarray*}
\label{eq:embed}
j: \mathcal{C}(\overline{K})  \rightarrow &  \Pico(\overline{K}) & = \J(\overline{K}) \\
    R           \mapsto & (R) - (Q)           
\end{eqnarray*}   
\\
where $Q \in \mathcal{C}(K)$ and $(R)$ denotes the divisor class of $R$. We note that since $Q \in \mathcal{C}(K)$ then if $R$ is a strictly-$L$ point on $\mathcal{C}$ then $j(R)$ is a strictly-$L$ point on $\J$. So all that remains is to construct the points on $\mathcal{C}$.\\
\\
It is a fact that every such $\mathcal{C}$ is birational over $K$ to a plane curve $\cp$ which has the same genus but is not necessarily smooth. We can resolve these singularities by blowing up to obtain a smooth curve $\cp_s$. Blowing up carries strictly-$L$ points on $\cp$ to strictly-$L$ points on $\cp_s$ and since $\mathcal{C}$ is birational to $\cp_s$ over $K$ we have $\rk(\mathcal{C} / L) = \rk (\cp_s/ L)$ for all $L/K$. Thus it is sufficient to construct strictly-$L$ points on $\cp$.\\
\\
We do this by constructing covering covers $\phi :\mathcal{C}_1 \rightarrow \cp$ where $\mathcal{C}_1$ has degree $p$ and where the map $\phi$ and the curves $\mathcal{C}_1$ are explicit. Indeed $\mathcal{C}_1$ and $\phi$ will be constructed so that the strictly-$L$ points we find on $\mathcal{C}_1(L)$ are in fact $S$-integers (for a fixed set of places $S$) and such that $\phi$ carries strictly-$L$ points on $\mathcal{C}_1$ to strictly-$L$ points on $\cp$. This will allow us to apply Siegel's Theorem to assert that $\mathcal{C}_1(L)$ can have only finitely many $S$-integer points and so we can deduce that the infinitely many points on $\mathcal{C}_1$ we generate must lie inside infinitely many different degree $p$ extensions, as desired.

\section{Proof of Theorem \protect ~\ref{th:jac}}

Let us suppose that we have an irreducible curve $\mathcal{C}$ defined over $K$ and with genus $g(\mathcal{C})\geq 1$. We shall suppose that $\mathcal{C}$ is birational to a plane curve $\cp$ of the form

\[ \cp: g(y) = f(x) \]
\\
where $f$ and $g$ are polynomials having degree $k$ and $d$ respectively. We suppose further that $(d,k)=1$ which allows us to assume that $f$ and $g$ are both monic and also, given some $q \in \OK_K$, that their non-leading coefficients are divisible by $q^n$ for any $n \geq 0$. Since $(d,k) = 1$ there are $a,b \in \Z$ such that $bk-ad=1$. Interchanging $x$ and $y$ if necessary, we can assume that $a,b >0$.\\
\\
Consider the following rational map:

\begin{eqnarray*}
\label{eq:phi2}
\phi : & \A^2_K & \rightarrow  \A^2_K \\
       &  (u,t) & \mapsto  \left(u + \frac{q^b}{t^{n}}, q^{a}t^{m} \right)
\end{eqnarray*} 
\\  
where $n,m>0$ are integers to be specified later and $q \in \OK_K$ is a generator of any prime ideal. We shall construct a cover of $\cp$ by taking the Zariski closure of the preimage of $\cp$ under $\phi$. Call this curve $\mathcal{C}_1$. It is cut out by the following equation:

\[ \mathcal{C}_1: h(u,t) :=  g(q^a t^m)t^{kn} - f(u)t^{nk} - f^{(1)}(u) q^{b}t^{n(k-1)} - \frac{f^{(2)}(u)}{2!} t^{n(k-2)}q^{2b} - \ldots - q^{kb} = 0. \]
\\
Note that $h$ is of degree $dm+kn$ in $t$, has leading coefficient $q^{ad}$ and has all its coefficients in $\OK_K$. The most important fact for us is the following.

\begin{lem}
\label{lem:irred}
h is irreducible over K.
\end{lem}

\proof{
Let us first of all assume (after perhaps performing a change of variables on our base curve $\cp$) that the non-leading coefficients of $g(y)$ and $f(x)$ are all divisible by $q^{da+1}$. Suppose that $h$ is reducible. Then for any specialisation of the variable $u$ to an element $m \in K$ the resulting polynomial $h(m,t)$ must be reducible into two polynomials in $t$ of degree at least one. Consider then the polynomial $h(0,t)$. Since $\frac{1}{i!}f^{(i)}(0)$ is just the $i^{th}$ coefficient of $f$ we see that every coefficient of $h(0,t)$ is divisible by $q^{da}$. Hence $h(0,t) = q^{da}h_1(t)$. Since $kb-ad=1$ we note that the constant term in $h_1(t)$ is $q$, the leading term is $1$ and every other coefficient is divisible by $q$. Hence $h_1(t)$ is Eisenstein and thus irreducible. Therefore, $h(0,t)$ and thus $h(u,t)$ are irreducible also. \\
\qed}
\\
We are now in a position to complete the proof of Theorem \protect ~\ref{th:jac}. \\
\\
$\textit{Proof of Theorem \protect ~\ref{th:jac}.}$ We want to find an $N(\mathcal{C}) > 0$ such that for all primes $p \geq N(\mathcal{C})$ we have an affirmative answer to Question \protect ~\ref{q:1}. We formed a cover $\mathcal{C}_1$ of $\cp$ which is given by $\mathcal{C}_1 : h(u,t) = 0$ where $h$ has degree $md+nk$. Now $(k,d) = 1$, therefore if we let $n$ range in $0 \leq n \leq d-1$ then we have that $nk$ occupies all the congruence classes modulo $d$. This is easily seen: if $ak \equiv bk (\text{mod }d)$ then $a \equiv b (\text{mod }d)$. Therefore we can express every number greater than $k(d-1)$ in the form $md+nk$ and in particular we can represent every prime $p > k(d-1)$ by $p = md + kn$ with $0 < m$ and $0 < n \leq d - 1$. \\
\\
Now we are in a good position. We can form a curve $\mathcal{C}_1: h(u,t) = 0$ of degree $p$ with $h(u,t) \in \OK_K$ and irreducible for all primes $p \geq k(d-1) + 1 = N(\mathcal{C})$. By Hilbert's Irreducibility Theorem we can find infinitely many $m \in K$ such that $h(m,t)$ is an irreducible polynomial of degree $p$ in $t$. Take such an $m$ and let $\alpha_m$ be a root of $h(m,t)$. Then $Q_m = (m,\alpha_m) \in \mathcal{C}_1(L_m)$ where $L_m = K(\alpha_m)$. Note that $[L_m:K] = p$. We can form infinitely many such points and since $\phi\vert_{\mathcal{C}_1}:\mathcal{C}_1 \rightarrow \cp$ has degree strictly less than $p$ we see that strictly-$L_m$ points on $\mathcal{C}_1$ are taken to strictly-$L_m$ points on $\cp$. (This is in virtue of the fact that if $p$ is prime and $l < p$ then $l$ and $p$ are coprime - this is the only place we use the primality of $p$ in our argument.)\\
\\
By Lemma \protect \protect ~\ref{lem:1}, all we need to show now is that the extensions $L_m$ do not coincide too often. To see this we note that if $S = \{\mathfrak{P} \in \text{Spec }\OK_{L_m} : q \in \mathfrak{P} \}$ then the roots of the polynomial $g(m,t)$ are all $S$-integers and by Siegel's Theorem there are only finitely many $S$-integers lying on $\cp$. Therefore, there must be infinitely many different extensions $L_m$ and so by Lemma \protect ~\ref{lem:1} we are done.
\qed\\
\\
We can now use this result to prove some interesting corollaries as stated in the introduction. 

\newtheorem{cor1}{Corollary}
\begin{cor1}
Let $\mathcal{C}$ be birational over $K$ to the plane curve cut out by the equation $\cp: y^2 = g(x)$ where $k=\text{deg}(g)$ is odd. Then there is an $N(\mathcal{C})$ such that for all primes $p \geq N(\mathcal{C})$ Question \protect ~\ref{q:1} has an affirmative answer. What is more, we can take $N(\mathcal{C}) = k+1$.  
\end{cor1}

\proof{
This is a straightforward application of Theorem \protect ~\ref{th:jac}. The proof of the theorem gives $N(\mathcal{C}) = k(d-1)+1 = k +1$.
\qed
}
This has the following pleasing corollary.

\begin{cor1}
Question 1 has an affirmative answer for every prime $p$ when $\mathcal{C}$ is an elliptic curve.
\end{cor1}

\proof{
Since elliptic curves can be written in the form $\mathcal{C} : y^2 = x^3 + Ax + B = g(x)$, Corollary \protect ~\ref{th:hyp} gives us the result for all primes $p \geq 5$. The primes 2 and 3 are not hard to prove. For $p=2$, if we take $m \in \OK_K$ then $P_m = (m, \sqrt{f(m)})$ lies on $\mathcal{C}$. $f(m)$ cannot be a square infinitely often or else $E(\OK_K)$ would be infinite, contradicting Siegel's Theorem. The same theorem also tells us that the points $P_m$ must lie in infinitely many different quadratic extensions. For $p=3$, we apply the same argument by putting $y = m \in \OK_K$. This gives a monic cubic polynomial in $x$. If this is reducible then there must be a linear factor and since the cubic is monic this yields an $\OK_K$ point on $\mathcal{C}$. Thus the cubics obtained this way must be irreducible almost all of the time and so give us our desired degree 3 points. \qed  \\
}
\\
$\textbf{Remark.}$ It is worth remarking that the degree $p$ extensions for $p>3$ which are constructed in the proof of Theorem \protect ~\ref{th:jac} take on, for elliptic curves, a particularly simple form. For a curve of the shape $y^2 = x^3 + Ax + B$ the constructed degree $p$ extensions are of the form $K[t]/g(k,t)$ where 

\[ g(k,t) = q^2t^p - (k^3+Ak+B)t^3 - (3k^2 + A)qt^2 - 3q^2kt - q^3, \]
\\
$q$ is any prime in $\OK_K$ and $k$ is to be considered as an element of $\OK_K$ so that $g(k,t) \in K[t]$. If we specify $k$ so that $g(k,t)$ is irreducible then we have that the extension is generated by a root of a polynomial of the form $q^2t^p - f(t)$ where $f \in \OK_K[t]$ is a cubic.

\section{A result when the degrees are not coprime}

A major hypothesis in the statement of Theorem \protect ~\ref{th:jac} is that the degrees of the two polynomials involved be coprime. In this final section we shall consider an example of a class of curves not fitting this restriction but where something can still be said. 

\begin{thm}
\label{th:hyp2}
Let $\mathcal{C}$ be a smooth, irreducible curve of positive genus, defined over $K$ and birational to a plane curve of the shape

\[ \cp : y^d = x^k + D \]
\\
where $D \in \OK_K^\times$ (and with no restrictions on $d$ or $k$). Then there is a $N(\mathcal{C}) > 0$ such that for all primes $p > N(\mathcal{C})$ there are infinitely many extensions $L/K$ of degree $p$ where $\rk(\mathcal{C}/L) > \rk(\mathcal{C}/K)$. 
\end{thm}

\proof{
We follow our usual strategy of constructing points on $\cp$. We begin by covering $\cp$ by the following curve

\[ \mathcal{C}_1 : y^n = x^n + D \]
\\
where $n = \text{lcm}(d,k)$. Note that we can assume that $q^n$ divides $D$ exactly for any $q$ outside of a finite set. Let us assume then that $q^{2n}$ exactly divides $D$ for some prime $q$. Let $w = y - x$, then by factorising $y^n-x^n$ we have

\[ \mathcal{C}_1 : w^n + w^{n-1}x + \ldots wx^{n-1} = D .\]
\\
Now, we are going to form a cover of $\mathcal{C}_1$ by looking at its preimage under the map $\phi : \A_{(u,t)}^2 \rightarrow \A_{(x,w)}^2$ given by

\[ \phi(u,t) = \left(qt^s + u , q^{n}t^r\right). \]
\\
As we have seen, the key point is showing that the curve $\mathcal{C}_2 : h(u,t) = 0$ given as the Zariski closure of $\phi^{-1}(\mathcal{C}_1)$ is irreducible. Consider the specialisation $h(0,t)$. This is given by

\[ h(0,t)  =   (q^nt^r)^n + (q^n t^{r})^{n-1}(qt^{s}) + \ldots + (q^nt^r)(qt^s)^{n-1} - D \]
\\
and one can observe that ever coefficient is divisible by $q^{2n}$ except for the coefficient of $t^{(n-1)s + r}$ (the last non constant monomial in the equation above). This has a coefficient of $q^{2n-1}$ and so we can divide out by this to get an Eisenstein polynomial, showing that $h(0,t)$ is irreducible and thus so is $h(u,t)$. \\
\\
If we have that $s > r$ in $\phi$ then $h(u,t)$ is of degree $(n-1)s + r$ in $t$. We can express every prime $p \geq (n-1)(n-2)$ as $p = (n-1)s + r$ for some $r, s$ with $s>r$ and so we see, by Hilbert's Irreducibility Theorem that for such a $p$ we can specialise $h(u,t)$ to $h(m,t)$ for infinitely many $m \in K$ to get an irreducible polynomial of degree $p$. The roots of this polynomial will then correspond to points strictly in a degree $p$ extension $L/K$. These points are then mapped to our original plane curve $\cp$. Since $p$ is larger than the degrees of any of the covering maps involved we see that the end point is still strictly-$L$. Thus we are done by Lemma \protect ~\ref{lem:1}. 
\qed}
\\


\begin{thebibliography}{9}

\bibitem{DD}
  Dokchitser, T., Dokchitser, V.,
  Elliptic curves with all quadratic twists of positive rank
  \emph{Acta Arithmetica},
  \textbf{139} (2009), 193-197.
\\
\bibitem{TD}
  Dokchitser, T.,
  Ranks of elliptic curves in cubic extensions,
  \emph{Acta Arithmetica},
  \textbf{126} (2007), 357-360. 
\\ 
\bibitem{FKK}
  Fearnley, J., Kisilevshy, H., Kuwata, M.,
  Vanishing and Non-Vanishing Dirichlet Twists of $L$-Functions of Elliptic Curves,
  \emph{Journal of the London Mathematical Society},
  first published online June 7, 2012 doi;10.1112/jlms/jds018.
\\
\bibitem{G}
  Goldfeld, D.,
  Conjectures on elliptic curves over quadratic fields,
  \emph{Springer Lect. Notes},
  \textbf{751} (1979), 108-118.   
\\  
    
\end{thebibliography}
\end{document}